\documentclass[ssy]{imsart}
\usepackage{amsmath,amsfonts,amsthm}
\usepackage{amssymb}

\usepackage{epsfig}

\usepackage{graphics}
\usepackage{graphicx}

\usepackage{latexsym}
\usepackage{graphics}
\usepackage{amsmath}
\usepackage{amsthm}
\usepackage{xspace}
\usepackage{amssymb}
\usepackage{color}
\usepackage{cite}
\usepackage{latexsym}
\usepackage{graphics}
\usepackage{amsmath}
\usepackage{amsthm}
\usepackage{xspace}
\usepackage{amssymb}
\usepackage{epsfig}

\newcommand{\beql}[1]{\begin{equation}\label{#1}}
\newcommand{\eeql}{\end{equation}}
\newcommand{\eqn}[1]{(\ref{#1})}

\newcommand{\R}{\mathbb{R}}

\newcommand{\E}{\mathbb{E}}

\newcommand{\cs}{{\cal S}}

\newcommand{\cx}{{\cal X}}

\newcommand{\Z}{\mathbb{Z}}

\newtheorem{thm}{Theorem}
\newtheorem{lem}[thm]{Lemma}

\newtheorem{cor}[thm]{Corollary}

\newtheorem{definition}[thm]{Definition}

\theoremstyle{remark}
\newtheorem{rem}[thm]{Remark}

\begin{document}

\begin{frontmatter}

\title{Tightness of stationary distributions of a flexible-server system
in the Halfin-Whitt asymptotic regime}
\runtitle{Tightness of stationary distributions}

\begin{aug}
  \author{\fnms{Alexander L.}  \snm{Stolyar}\ead[label=e1]{stolyar@lehigh.edu}}

  \runauthor{A. L. Stolyar}

  \affiliation{Lehigh University}

  \address{Bethlehem, PA\\
          \printead{e1}}

\end{aug}

\begin{abstract}

We consider a large-scale flexible service system 
with two large server pools and two types of customers.
Servers in pool 1 can only serve type 1 customers, 
while server in pool 2 are flexible -- they can serve both types 1 and 2.
(This is a so-called ``N-system.'' Our results hold for a more general
class of systems as well.) 
The service rate of a customer depends both on its type and the pool where it is served.
We study a priority service discipline, where 
type 2 has priority in pool 2, and type 1 prefers pool 1.
We consider the Halfin-Whitt asymptotic regime, where the arrival rate 
of customers and the number of servers in each pool increase to infinity 
in proportion to a scaling parameter $n$, while the overall system 
capacity exceeds its load by $O(\sqrt{n})$.

For this system we prove tightness of diffusion-scaled 
stationary distributions. Our approach relies on a single common
Lyapunov function $G^{(n)}(x)$, depending on parameter $n$ and defined on the entire state space
as a functional of the {\em drift-based fluid limits} (DFL).
Specifically, $G^{(n)}(x)=\int_0^\infty g(y^{(n)}(t)) dt$, where $y^{(n)}(\cdot)$ is the DFL
starting at $x$, and $g(\cdot)$ is a ``distance'' to the origin. ($g(\cdot)$ is same for all $n$).
The key part of the analysis is the study of the (first and second) derivatives of the DFLs and function $G^{(n)}(x)$.
The approach, as well as many parts of the analysis, are quite generic and
may be of independent interest.

\end{abstract}



\end{frontmatter}

\maketitle


\section{Introduction}

In this paper we consider a large-scale service system in the 
so-called Halfin-Whitt asymptotic regime. Such systems
received a lot of attention in the literature, especially in the past 10-15 year,
because they find a variety of applications, including, e.g., large customer contact centers 
\cite{AksinArmonyMehrotra,GansKooleMandelbaum} and large computer farms in network clouds.
The Halfin-Whitt regime, introduced originally in \cite{HalfinWhitt81},
is such that the system capacity (roughly, number of servers) increases 
in proportion to a scaling parameter $n$, and exceeds the system load 
by $O(\sqrt{n})$. It is attractive because it allows -- in principle, under a good control 
algorithm -- to achieve both good performance (e.g. waiting times) and high
resource utilization.

In the Halfin-Whitt regime, the stochastic process describing the system behavior
is usually studied under {\em diffusion scaling}, i.e. it is 
centered at the system equilibrium point and scaled down by $n^{-1/2}$.
This name reflects the fact that, in the limit on $n\to\infty$, 
on any finite time interval, the sequence of diffusion-scaled processes $Y^{(n)}(\cdot)$
``typically'' converges to a diffusion process $Y(\cdot)$. Then, a fundamental 
question is whether or not the following 
{\em limit interchange} property holds: {\em the limit of stationary distributions
of $Y^{(n)}(\cdot)$ is equal to the stationary distribution of $Y(\cdot)$.}
In turn, the key difficulty in establishing
 the limit interchange property, is verifying
the {\em (stationary distribution) tightness} property:
 {\em the family of stationary distributions of $Y^{(n)}(\cdot)$ is tight}. 

The tightness property in the Halfin-Whitt regime is usually difficult to verify even for 
systems with single pool of homogeneous servers, if there is more than
one type of arriving customers and/or the service time distribution is 
 non-exponential; 
see \cite{Gamarnik_Momcilovic, GamarnikGoldbergGGN, Gamarnik_Stolyar, DaiDiekerGao-tightness} for the results in this direction.
(We note that the problems of verifying
the tightness and limit interchange
exists not only in the Halfin-Whitt regime, but also in the so-called
{\em conventional heavy traffic} regime; 
see \cite{Gamarnik_Zeevi, BudhirajaLee, Gurvich-HT}.)
More general models, where 
there are multiple flexible server pools with different capabilities (service rates)
w.r.t. different customer types, pose additional challenges.
The key additional difficulty is that for such systems the state space is ``fractured''
into multiple domains, where the process dynamics is very different.
Papers \cite{SY2010, SY2012, St2013_tightness} contain tightness /
limit interchange results for some flexible multi-pool models; 
although, \cite{SY2012, St2013_tightness} consider a strictly subcritical load
regime (different from Halfin-Whitt), in which the capacity exceeds
the load by $O(n)$.

One approach for verifying the stationary distribution tightness
is to find a single common Lyapunov function, for which an appropriate
``negative expected drift'' condition can be established.
This approach is used in 
\cite{Gamarnik_Momcilovic, Gamarnik_Stolyar, SY2010, DaiDiekerGao-tightness}.
(Papers \cite{GamarnikGoldbergGGN, SY2012, St2013_tightness} use different approaches,
not relying on a single Lyapunov function.)
Of course, finding/constructing a suitable Lyapunov function is usually the key challenge.
For example, paper \cite{DaiDiekerGao-tightness}, which proves tightness
for a single-pool model with first-come-first-serve 
discipline and
phase-type service time distribution, uses an elaborate
{\em common quadratic Lyapunov function}, of the type proposed in 
\cite{DiekerGao-piecewise}; the tightness result in \cite{DaiDiekerGao-tightness} 
also requires that the customers waiting in the queue abandon at positive rate.
And again, finding single common Lyapunov function is further complicated for flexible multi-pool systems.

\subsection{Paper contributions}

We consider a flexible multi-pool system
with two customer types and two server pools 
(the so-called
``$N$-system''), under a priority discipline, in the Halfin-Whitt regime,
and prove the stationary distribution
tightness result, Theorem~\ref{th-tightnessN}, 
which implies the limit-interchange, Corollary~\ref{cor-interchange}.
(These results hold for a more general
class of systems as well, as discussed in Section~\ref{sec-generalization}.) 

The state space of the diffusion-scaled process for $N$-system
has five domains, where the process drift is given by different affine functions;
the domain boundaries depend on parameter $n$.
Nevertheless, 
we construct a single Lyapunov function $G^{(n)}(x)$ (depending on parameter $n$) on the entire state space,
as a functional of the {\em drift-based fluid limits} (DFL),
which are the deterministic trajectories defined by the drift of the process.
Specifically, 
\beql{eq-Gdef-contrib}
G^{(n)}(x)=\int_0^\infty g(y^{(n)}(t)) dt, 
\end{equation}
where $y^{(n)}(\cdot)$ is the DFL
starting at $x$, and $g(\cdot)$ is a ``distance'' to the origin. (Function $g(\cdot)$ does {\em not} depend on $n$.)
For a Lyapunov function of this type, 
in a setting more general than needed for the proof of Theorem~\ref{th-tightnessN},
we give sufficient conditions for the 
tightness in Theorem~\ref{th-criterion}; the key condition a (uniform in $n$) bound 
on the Lyapunov function second derivatives. This result may be of independent interest.

The proof of Theorem~\ref{th-tightnessN} verifies the 
conditions of Theorem~\ref{th-criterion} for the $N$-system.
This requires the analysis of the DFL structure, and of the (first and second) derivatives of DFLs 
and corresponding functionals $G^{(n)}(x)$ 
on the initial state $x$; it also requires an appropriate choice of the "distance" $g(\cdot)$. Many parts of this analysis are quite generic
and may also be of independent interest.

As will be illustrated below in Section~\ref{sec-intuition},
for a {\em deterministic} dynamic system, with trajectories $y^{(n)}(\cdot)$ defined
by a continuous derivative-field,
the function $G^{(n)}(x)$ given by \eqn{eq-Gdef-contrib}
is a natural Lyapunov function,
as long as it is well defined (the integral in \eqn{eq-Gdef-contrib} is finite).
In queueing networks literature, this observation is used, for example, in \cite{Ye_Chen_2001, Schonlein_Wirth} to establish the existence 
of a Lyapunov function for stable {\em deterministic} fluid models.
This observation, however, does {\em not} imply that $G^{(n)}(x)$ defined by \eqn{eq-Gdef-contrib} via DFLs $y^{(n)}(\cdot)$
can serve as a Lyapunov function 
for a (family of) {\em random} process(es). 
In this paper we give sufficient conditions under which Lyapunov functions $G^{(n)}(x)$ can be used 
to establish tightness of stationary distributions, and then verify these conditions for the $N$-system.

A Lyapunov function similar in spirit to \eqn{eq-Gdef-contrib} was used in \cite{DW94} to establish a sufficient 
condition for positive recurrence of a semimartingale reflecting Brownian motion in the positive orthant.
(In \cite{DW94}, the solutions to the Skorohod problem, for the trajectories determined by the process drift alone, are the DFLs in our terminology.)
Obtaining the Lyapunov function second derivative bounds is also a key part of the analysis in \cite{DW94}.
We note, however, that our basic model, the problem, the structure of the (family of) process(es) and corresponding DFLs, the form of function $g(\cdot)$, and the analysis
of the Lyapunov function derivatives are completely different.

\subsection{Layout of the rest of the paper} 

In Section~\ref{sec-intuition}, we informally discuss our 
general approach and the Lyapunov function construction.
Section~\ref{sec-N-ssytem-def} formally defines the $N$-system, the Halfin-Whitt
regime for it, and states the tightness (Theorem~\ref{th-tightnessN}) and
the limit-interchange (Corollary~\ref{cor-interchange}) results.
In Section~\ref{sec-tightness-criterion}, in a setting more 
general than needed for the $N$-system,
we give a formal construction of the DFLs
and the Lyapunov function, and sufficient conditions for the 
tightness (Theorem~\ref{th-criterion}). 
Section~\ref{sec-proof-N} contains the proof of Theorem~\ref{th-tightnessN};
here we choose a specific ``distance'' function $g$ and verify the 
conditions of Theorem~\ref{th-criterion} for the $N$-system.
A generalization of the $N$-system, for which our results still hold,
is described in Section~\ref{sec-generalization}.
Finally, in Section~\ref{discussion}, we discuss our approach and results.

\subsection{Basic notation}

Symbols $\R, \R_+, \Z, \Z_+$ denote the sets of real, real non-negative,
integer, and integer non-negative numbers, respectively.
In the Euclidean space $\mathbb{R}^I$ (of dimension $I\ge 1$):
$|x|$ denotes standard Euclidean
norm of vector $x=(x_1,\ldots,x_I)$, while 
$\|x\| = \sum_i |x_i|$ denotes its $L_1$-norm;
scalar product of two vectors is denoted $x \cdot y = \sum_i x_i y_i$;
$diag (x)$ denotes diagonal square matrix, with diagonal elements given by $x$;
we write simply $0$ for a zero matrix or vector;
vectors are written as row-vectors, but in matrix expressions they are viewed
as column-vectors (without using a transposition sign). 
For real numbers $u$ and $w$: 
$u \vee w = \max\{u,w\}$, $u \wedge w = \min\{u,w\}$, and $\lfloor u \rfloor$ denotes the largest integer not greater than $u$. 

For a vector-function $ y(\cdot)=(y(t), ~t\ge 0)$, we denote
$\|y(\cdot)\| = \sup_{[0,\infty)} \|y(t)\|$. Abbreviation {\em u.o.c.} means
{\em uniform on compact sets} convergence. 
If $X(t), ~t\ge 0,$ is a Markov process, 
we write $X(\infty)$ for a random element with the distribution equal to a
stationary distribution of the process. (In all cases considered in this paper, the stationary distribution will be unique.)
Symbol $\Rightarrow$ denotes convergence in distribution of random elements;
random processes are random elements in the appropriate Skorohod space.
For a condition/event $H$, 
the indicator function $I\{H\}$ is equal to $1$ when $H$ holds and $0$ otherwise.

\section{The intuition for the Lyapunov function construction}
\label{sec-intuition}

The discussion in this entire section is informal.
Consider a deterministic dynamic system governed by ODE 
\beql{eq-diff777}
(d/dt)y = v(y),
\end{equation}
where state $y$ is a vector, and the vector-field $v(\cdot)$ is Lipschitz continuous. 
Suppose the system has unique stable point $0$.
Let $g(x)$ be a non-negative continuous (and sufficiently smooth) function, 
which measures a "distance" from $0$.
(In our results, we will use $g(x)$ which is a smooth approximation of 
$L_1$-norm $\|x\|$.)
Suppose that for any initial state $y(0)=x$ the
trajectory $y(t), ~t\ge 0$ converges to $0$ and, moreover, 
\beql{eq-lyap-def}
G(x) = \int_0^\infty g(y(t)) dt < \infty.
\end{equation}
Then $G(\cdot)$ is a Lyapunov function for this dynamic system, in the sense that
$$
(d/dt) G(y(t)) = G'(y(t)) \cdot v(y(t)) = -  g(y(t)),
$$
where $G'$ denotes the gradient of $G$. (The first equality is immediate from \eqn{eq-diff777}, and $(d/dt) G(y(t)) = -  g(y(t))$ is from
$G(y(t)) = \int_t^\infty g(y(s)) ds$.)

Suppose now that instead of a deterministic system we have a Markov process $Y(\cdot)$,
 for which vector-field $v(\cdot)$ gives the drift.
Then we can define deterministic trajectories $y(\cdot)$,
 and function $G(\cdot)$, the same way as above.
(The trajectories $y(\cdot)$ we call drift-based fluid limits (DFL).)
Suppose further that the process generator $A$ is such that
\beql{eq-generator-intuition}
A G(y) = G'(y) \cdot v(y) + H(y),  ~~~~ |H(y)| \le C_0 \|G''(y)\|,
\end{equation}
where $C_0$ is a constant and $G''$ denotes the Hessian matrix of second derivatives.
(To interpret \eqn{eq-generator-intuition} one can think, for example,
of a diffusion process with bounded diffusion coefficients.
In this paper we will work {\em not} with diffusion processes, 
but rather with diffusion-scaled processes for our queueing system -- 
their behavior can be very different
from that of diffusions, especially when the system state is "far" from 
the equilibrium point. Nevertheless, the process generator will have 
form \eqn{eq-generator-intuition}.) Then, we have
$$
A G(y) \le G'(y) \cdot v(y) + C_0 \|G''(y)\| = - g(y) + C_0 \|G''(y)\|.
$$
If we can show that 
\beql{eq-second-deriv-intuition}
\|G''(y)\| \le C_1 g(y) + C_2
\end{equation}
with a sufficiently small $C_1$, then for some $\epsilon>0$ and $\kappa>0$,
\beql{eq-second-deriv-intuition2}
A G(y) \le -\epsilon g(y) + \kappa.
\end{equation}
This is a Lyapunov-Foster type condition from which we can obtain the steady-state bound
$\E g(Y(\infty)) \le \kappa/\epsilon$, where $Y(\infty)$ is $Y(t)$ when the process is in stationary regime.

Finally, suppose we consider a family of processes $Y(\cdot)=Y^{(n)}(\cdot)$,
with the drift $v(\cdot)$ and
generator $A$ depending on $n$. If for some 
common function $g$ such that $g(x)\to \infty$ as $\|x\|\to \infty$,
we can derive estimates \eqn{eq-generator-intuition}-\eqn{eq-second-deriv-intuition2}
with constants independent of $n$, then $\E g(Y^{(n)}(\infty))$ is bounded uniformly in $n$,
and therefore the family of stationary distributions of $Y^{(n)}(\cdot)$ is tight.

This is the program that we implement in this paper, for the sequence of
diffusion-scaled processes for the $N$-system.
The difficult part is obtaining the second derivative bound 
\eqn{eq-second-deriv-intuition}.
Since $G$ is defined as a functional of the DFLs $y(\cdot)$,
this involves the analysis of the dependence of DFLs on the initial state.

\section{$N$-system with absolute priority}
\label{sec-N-ssytem-def}

\begin{figure}[htbp]
\begin{center}
\includegraphics[width=90mm]{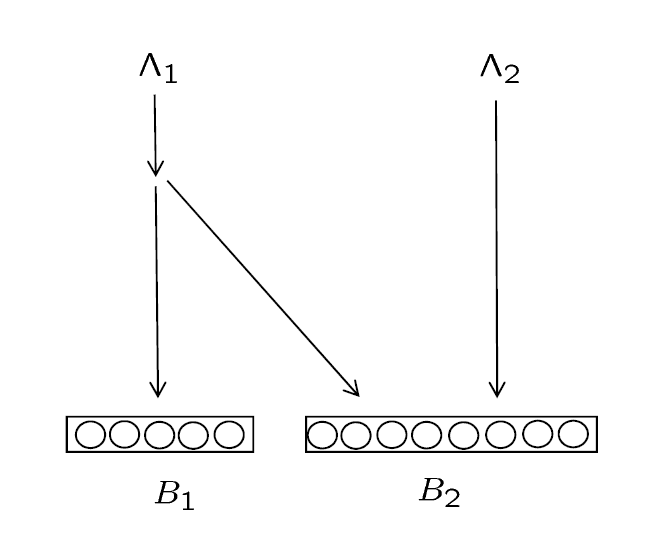}
\end{center}
\vspace{-0.2in}
\caption{$N$-system.}
\label{fig:N-system}
\end{figure}

Consider a so-called $N$-system, with absolute priorities. 
(See Fig.~\ref{fig:N-system}.)
There are two customer types, arriving according to as independent Poisson 
processes with rates $\Lambda_1>0$ and $\Lambda_2>0$, respectively.
There are two server pools, with $B_1$ and $B_2$ identical servers, respectively.
The total service requirement of any customer is an independent, exponentially distributed random variable with mean $1$.
A customer of type 2 can only be served by a server in pool 2, and if it does receive service, it does so at rate $\mu_{22}>0$. 
A customers of type 1 can be served by a server in either pool 1 or 2, with service rates being $\mu_{11}>0$ and $\mu_{12}>0$,
respectively. Type 2 customers have absolute (preemptive) priority (in pool 2); namely, if there are $X_2$ type 2 customers in the
system, as many of them as possible, $X_2 \wedge B_2$,  receive service in pool 2, and the remaining $X_2 - X_2 \wedge B_2 = (X_2-B_2) \vee 0$
wait in the queue. (Here $\wedge$ and $\vee$ denote minimum and maximum, respectively.) 
Therefore, the total
service rate of all type 2 customers is 
\beql{eq-serv-rate2}
\mu_{22} (X_2 \wedge B_2).
\end{equation}
The type 1 customers have absolute preference
to be served in pool 1, and have lower preempt-resume priority in pool 2. Namely, if there are $X_1$ type 1 customers in the system,
then  $X_1 \wedge B_1$ of them are served in pool 1, $[(X_1-B_1) \vee 0] \wedge [(B_2-X_2) \vee 0]$ are served in pool 2,
and the remaining $[X_1 - (B_1+B_2) + (X_2 \wedge B_2)] \vee 0$ wait in queue. The total
service rate of all type 1 customers is 
\beql{eq-serv-rate1}
\mu_{11} \{X_1 \wedge B_1\} + \mu_{12} \{[(X_1-B_1) \vee 0] \wedge [(B_2-X_2) \vee 0]\}.
\end{equation}

\begin{figure}[htbp]
\begin{center}
\includegraphics[width=90mm]{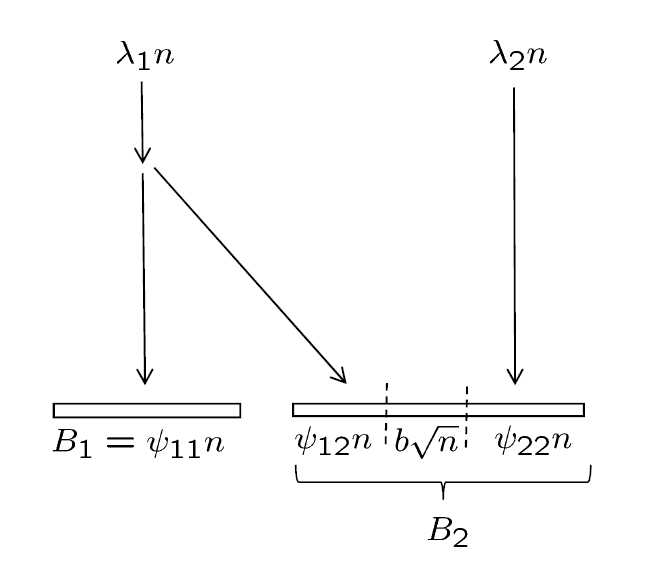}
\end{center}
\vspace{-0.2in}
\caption{$N$-system in Halfin-Whitt asymptotic regime.}
\label{fig:N-system-sequence}
\end{figure}

We  consider a sequence of such systems, indexed by
a positive scaling parameter $n$, increasing to infinity. 
(See Fig.~\ref{fig:N-system-sequence}.)
In a system with parameter $n$,
\beql{eq-param1}
\Lambda_1 = \lambda_1 n, ~~ \Lambda_2 = \lambda_2 n, 
\end{equation}
\beql{eq-param2}
 B_1= \psi_{11} n, ~~ B_2= \psi_{12} n + \psi_{22} n + b \sqrt{n},
\end{equation}
where the positive parameters $b, \lambda_1, \lambda_2, \psi_{11}, \psi_{12}, \psi_{22}$ are such that
\beql{eq-param3}
 \lambda_2  = \psi_{22}  \mu_{22}, ~~ \lambda_1 = \psi_{11}  \mu_{11} + \psi_{12}  \mu_{12}.
\end{equation}

Given this definition, and the priorities, the system
``desired operating point,'' which we will refer to as {\em equilibrium point},
is such that $X_2=\psi_{22} n$
and $X_1=\psi_{11} n + \psi_{12} n$, where type 1 customers occupy the entire pool 1
and $\psi_{12} n$ servers in pool 2; the equilibrium point is such that 
$b\sqrt{n}$ servers in pool 2 are idle -- this is the ``margin'' by which system
capacity exceeds its load. 
(Again, see Fig.~\ref{fig:N-system-sequence}.)

\begin{rem}
\label{remark-rounding}
To be precise, in the definition of the sequence of systems, we need to make sure that $B_1$ and $B_2$ are integer.
Equations \eqn{eq-param2}, as written above, assume that $B_1$ and $B_2$ "happen to be" integer.
We make this assumption throughout the paper to simplify the exposition, while maintaining rigor of the results and arguments.
More specifically, we could replace \eqn{eq-param2} with, for example,
\beql{eq-param222}
 B_1= \lfloor \psi_{11} n \rfloor, ~~ B_2= \lfloor\psi_{12} n + \psi_{22} n + b \sqrt{n} \rfloor.
\end{equation}
If we do that, it is easy to check that for each $n$ we can choose numbers $\psi_{ij}^{(n)}$, $(ij)=(11), (12), (22)$, 
and $b^{(n)}$, such that:  $|\psi_{ij}^{(n)} - \psi_{ij}| \le \kappa /n$ and $|b^{(n)}-b| \le \kappa/\sqrt{n}$
for some constant $\kappa>0$;
\eqn{eq-param222} can be rewritten as
\beql{eq-param222777}
 B_1= \psi_{11}^{(n)} n, ~~ B_2= \psi_{12}^{(n)} n + \psi_{22}^{(n)} n + b^{(n)} \sqrt{n};
\end{equation}
and \eqn{eq-param3} can be rewritten as 
\beql{eq-param333}
 \lambda_2  = \psi_{22}^{(n)}  \mu_{22}, ~~ \lambda_1 = \psi_{11}^{(n)}  \mu_{11} + \psi_{12}^{(n)}  \mu_{12}.
\end{equation}
The sequence of systems will then be defined by \eqn{eq-param1}, \eqn{eq-param222777}, \eqn{eq-param333}. 
Then, the entire analysis in this paper will hold as is, with $\psi_{ij}$ and $b$ replaced 
everywhere with $\psi_{ij}^{(n)}$ and $b^{(n)}$, respectively. (We note that
the components of the equilibrium point, namely 
$X_2=\psi_{22}^{(n)} n$
and $X_1=\psi_{11}^{(n)} n + \psi_{12}^{(n)} n$, need {\em not} be integer.)
\end{rem}

It is easy to see that for each $n$ the process $X^{(n)}(t)= (X^{(n)}_1(t),X^{(n)}_2(t)), ~t\ge 0$, is continuous-time countable irreducible Markov chain,
with the state space being (for each $n$) $\Z_+^2$. Further,
it is not difficult to check 
that, for each sufficiently large $n$, 
this Markov process is positive recurrent, and therefore has unique stationary distribution. Indeed, due to absolute priority, type 2 customers ``do not see'' 
type 1, and therefore $X_2^{(n)}(\cdot)$ in itself is a positive recurrent Markov chain,
which in steady-state occupies on average $\psi_{22} n$ servers in pool 2.
This means that on average $\psi_{12} n + b\sqrt{n}$ servers in pool 2
are available to serve type 1 customers; this is in addition to all
$\psi_{11} n$ servers in pool 1 which are available exclusively to type 1;
therefore, the average total service capacity available to type 1 is
$$
\psi_{11} n \mu_{11} + (\psi_{12} n + b\sqrt{n})\mu_{12} 
= \lambda_1 n + b \mu_{12} \sqrt{n} > \lambda_1 n.
$$
More details of a positive recurrence proof are given in Appendix \ref{app1}.

The diffusion-scaled version $\hat X^{(n)}(t)= (\hat X^{(n)}_1(t), \hat X^{(n)}_2(t))$ of the process $X^{(n)}(t)$ 
is defined by centering at the equilibrium point and rescaling by $1/\sqrt{n}$:
\beql{eq-diff-scale}
\hat X_1^{(n)} =(X_1^{(n)} -  \psi_{11}n -\psi_{12}n)/\sqrt{n}, ~~ \hat X_2^{(n)} =(X_2^{(n)} -  \psi_{22}n)/\sqrt{n}.
\end{equation}

\begin{thm}
\label{th-tightnessN}
For some $C>0$ and all sufficiently large $n$,
$$
\E  \|\hat X^{(n)}(\infty)\| \le C.
$$
\end{thm}

The proof of Theorem~\ref{th-tightnessN} is given in the rest of this paper.
It relies on a Lyapunov function (depending on $n$),
being a functional of a fluid trajectory,
determined by the process drift. Such fluid trajectories will be referred to as
drift-based fluid limits (DFL). In the rest of this section we define 
DFLs for the $N$-system under consideration, and give motivation 
for the form of Lyapunov function. Then, in Section~\ref{sec-tightness-criterion},
we give the Lyapunov function construction and sufficient tightness
conditions (Theorem~\ref{th-criterion}) in a setting that is more general 
than needed for the $N$-system. In the following sections we verify
the conditions of Theorem~\ref{th-criterion} for the $N$-system,
thus proving Theorem~\ref{th-tightnessN}.

For each $n$, for the unscaled process $X^{(n)}(\cdot)$, 
we define a drift function (vector field) $V^{(n)}=(V_1^{(n)},V_2^{(n)})$ for $x=(x_1,x_2) \in \R_+^2$.
(Note that it is defined on $\R_+^2$, and not just on the lattice $\Z_+^2$.) It is defined in the natural way, as the difference of arrival and service rates 
(see \eqn{eq-serv-rate2}-\eqn{eq-serv-rate1}):
\beql{eq-Vdef1}
V_1^{(n)} = V_1^{(n)}(x) = \Lambda_1  - \mu_{11} \{x_1 \wedge B_1\} - \mu_{12} \{[(x_1-B_1) \vee 0] \wedge [(B_2-x_2) \vee 0]\},
\end{equation}
\beql{eq-Vdef2}
V_2^{(n)} = V_2^{(n)}(x) = \Lambda_2 - \mu_{22} (x_2 \wedge B_2),
\end{equation}
where $\Lambda_1, \Lambda_2, B_1, B_2$ are the functions of $n$ given in \eqn{eq-param1}-\eqn{eq-param3}.

Let us denote by  $L_n$ the affine mapping $X^{(n)} \to \hat X^{(n)}$,
defined by \eqn{eq-diff-scale}.
Then, the state space of $\hat X^{(n)}$ is $\cs^{(n)} \equiv L_n \Z_+^2 \subset \cx^{(n)} \equiv L_n \R_+^2 \subset \R^2$.
Specifically, $\cx^{(n)}=\{x ~|~ x_1 \ge -  \psi_{11}\sqrt{n} -\psi_{12}\sqrt{n}, ~
x_2 \ge  -\psi_{22}\sqrt{n}\}$. The drift function for $\hat X^{(n)}$ is defined accordingly:
$$
v^{(n)}(x) = (1/\sqrt{n}) V^{(n)}(L_n^{-1} x), ~~ x \in \cx^{(n)}.
$$
We emphasize that $v^{(n)}(x)$ is defined on the continuous convex set $\cx^{(n)}$, which contains the discrete state space $\cs^{(n)}$.
It is important, however, that {\em at each point $x \in \cs^{(n)}$, $v^{(n)}(x)$ gives exactly the average drift of the process.}
Namely, 
\beql{eq-exact-drift}
v^{(n)}(x) = \sum_{x'} (x'-x) \nu^{(n)}(x,x'),
\end{equation}
where $\nu^{(n)}(x,x')$ is the Markov process transition rate from state $x$ to state $x'$; note that there is only a finite number
of "neighbor" states $x'$ for which $\nu^{(n)}(x,x')>0$.

As $n\to\infty$, set $\cx^{(n)}$ monotonically increases and converges to $\R^2$.

It is easy to observe that $v^{(n)}(x) = 0$ if and only if $x=0$; also, uniformly in $n$, $v^{(n)}(x)$ is Lipschitz continuous.
Given Lipschitz continuity of $v^{(n)}$, for any $x\in \cx^{(n)}$ there is a unique solution $y^{(n)}(t), ~t\ge 0,$ to the ODE
$$
(d/dt) y^{(n)}(t) = v^{(n)}(y^{(n)}(t)), ~~ y^{(n)}(0)=x.
$$
The solution stays within $\cx^{(n)}$ for all $t\ge 0$. (Indeed, for each $i$, $v_i^{(n)}(x) = \lambda_i \sqrt{n}$ when $x_i$
is at its lower bound -- see the definition of 
$\cx^{(n)}$.)
This trajectory $y^{(n)}(t), ~t\ge 0$, will be called the
{\em drift-based fluid limit} (DFL), starting from $x$. 

As we will show later in Section~\ref{sec-basic-dfl}, each DFL is such that $y^{(n)}(t)\to 0$ 
as $t\to\infty$. Moreover, after a finite time this convergence is 
exponentially fast, so that
$$
\int_0^\infty \|y^{(n)}(t)\| dt < \infty.
$$
The Lyapunov function we will use to prove Theorem~\ref{th-tightnessN} is
$$
G^{(n)}(x) = \int_0^\infty g(y^{(n)}(t)) dt < \infty,
$$
where $y^{(n)}(\cdot)$ is the DFL starting from $x$, and
$g(\cdot)$ is a smooth non-negative function (common for all $n$) 
approximating $\|\cdot\|$.

\begin{rem}
\label{remark-DFL}
In the literature on the steady-state tightness in the Halfin-Whitt asymptotic regime, 
deterministic trajectories defined by the drift vector field, which we call DFLs,
are considered in e.g. \cite{DaiDiekerGao-tightness},
where they are called fluid models. However, the way we use DFLs in this paper --
namely, to directly construct a Lyapunov function from them -- is completely different
from their use in \cite{DaiDiekerGao-tightness}.
\end{rem}

\subsection{Limit interchange}
\label{sec-interchange}

We conclude this section by noting that 
the tightness of stationary distributions 
of the processes $\hat X^{(n)}(\cdot)$, which follows from
Theorem~\ref{th-tightnessN}, allows us to easily establish
the limit interchange result, given in Corollary~\ref{cor-interchange} below. 

Recall that for a given $n$, the drift function $v^{(n)}(\cdot)$ is defined on the set 
$\cx^{(n)}=\{x ~|~ x_1 \ge -  \psi_{11}\sqrt{n} -\psi_{12}\sqrt{n}, ~ x_2 \ge  -\psi_{22}\sqrt{n}\}$,
which monotonically increases and converges to $\R^2$. It is easy to observe that $v^{(n)}(x) \to v(x)$ uniformly on compact subsets of 
$\R^2$,  where $v_2(x)=-\mu_{22} x_2$
and
$
v_1(x) = -\mu_{12} [x_1 \wedge (b-x_2)].
$
In fact, even stronger property holds: on any compact subset of 
$\R^2$, $v^{(n)}(x) = v(x)$ for all large $n$. 

\begin{cor}
\label{cor-interchange}
The following convergence holds
\beql{eq-lim-interchange}
\hat X^{(n)}(\infty) \Rightarrow \hat X(\infty),
\end{equation}
where $\hat X(\cdot)$ is a diffusion process which is a strong solution of SDE
\beql{eq-SDE}
d(\hat X) = v(\hat X)dt + (\sigma_1 dW_1,\sigma_2 dW_2),
\end{equation}
where $W_1, W_2$ are independent standard Brownian motions and
the diffusion coefficients are
$\sigma_1=[\lambda_1 + \psi_{11}\mu_{11}+ \psi_{12}\mu_{12}]^{1/2}$,
$\sigma_2=[\lambda_2 + \psi_{22}\mu_{22}]^{1/2}$.
\end{cor}

The proof is fairly straightforward, we just give an outline.
First, the following convergence on a finite interval holds
(see e.g. \cite{Gurvich_Whitt}).  Namely, 
consider a sequence of processes $\hat X^{(n)}(\cdot)$ with fixed 
initial states $\hat X^{(n)}(0) \to x\in \R^2$. Then, for any fixed $T_0>0$
\beql{eq-conv-to-diff}
(\hat X^{(n)}(t), ~t\in [0,T_0]) \Rightarrow (\hat X(t), ~t\in [0,T_0]),
\end{equation}
where $\hat X(\cdot)$ is a strong solution of \eqn{eq-SDE} with 
initial state $\hat X(0)= x$. Then, \eqn{eq-lim-interchange} can be
established, together with the existence and uniqueness of a stationary distribution 
of $\hat X(\cdot)$, as follows. We consider the sequence of stationary versions 
of the processes $\hat X^{(n)}(\cdot)$ on a fixed finite time interval
$[0,T_0]$, and let $n\to\infty$. Given tightness of stationary distributions
of pre-limit processes, we can choose a subsequence along which
$\hat X^{(n)}(0) \Rightarrow \tilde X(0)$ for some random vector $\tilde X(0)$;
then we also have $\hat X^{(n)}(T_0) \Rightarrow \tilde X(0)$. We then 
use \eqn{eq-conv-to-diff} to show that the distribution of $\tilde X(0)$
must be a stationary distribution of $\hat X(\cdot)$. The uniqueness
of the latter stationary distribution is easy to establish, for example, 
using a coupling argument.

\section{Lyapunov function construction and a tightness criterion}
\label{sec-tightness-criterion}

The model in this section is quite general (including the N-system as a special case). 
For this model we define DFLs, construct a functional of DFL, and give sufficient conditions under which this functional can serve 
as a Lyapunov function to prove tightness of stationary distributions. 
The section is self-contained, because its main construction and result may be of independent interest.
However, it may help the reader to keep the N-system described in 
Section~\ref{sec-N-ssytem-def} in mind as an example, to make the material
more concrete.

\subsection{Setting and assumptions}

Let $I \ge 1$ be a fixed positive integer. For each positive integer $n\ge n_0$ (where number $n_0$ is fixed), 
we consider a Markov chain $\hat X^{(n)}(t), ~t\ge 0$, with a countable state space
  $\cs^{(n)}$ which has the form
$$
\cs^{(n)} = \{ L_n x ~|~ x \in \Z^I\} \cap \cx^{(n)},
$$
where $\cx^{(n)}$ is a convex closed subset of $\R^I$, containing $0$, and $L_n x = x/\sqrt{n} + s^{(n)}$
with some fixed $s^{(n)}\in \R^I$. Assume that for each $n$ this Markov chain is irreducible, positive recurrent, and is such that the total transition rate out of any state is upper 
bounded by $R_1 n$ and any single transition has the jump size of at most $R_2/\sqrt{n}$,
where $R_1, R_2$ are positive constants independent of $n$. 
Suppose that, defined on $\cx^{(n)}$ is a drift function (vector field) $v^{(n)}(x)$, which is 
Lipschitz continuous uniformly in $n$. Assume that
 {\em at each point $x \in \cs^{(n)}$, $v^{(n)}(x)$ gives exactly the average drift of the process.}
Namely, 
\beql{eq-exact-drift-gen}
v^{(n)}(x) = \sum_{x'} (x'-x) \nu^{(n)}(x,x'),
\end{equation}
where $\nu^{(n)}(x,x')$ is the Markov process transition rate from state $x$ to state $x'$; given the upper bound on a single jump size,
 note that there is only a finite number
of "neighbor" states $x'$ for which $\nu^{(n)}(x,x')>0$.

Assume that for any $x\in \cx^{(n)}$, there is a unique solution  $y^{(n)}(t), ~ t\ge 0,$ to the ODE
$$
(d/dt) y^{(n)}(t) = v^{(n)}(y^{(n)}(t)), ~~ y^{(n)}(0)=x,
$$
and the solution stays within $\cx^{(n)}$. This solution is called the {\em drift-based fluid limit} (DFL), starting from $x$. 

Suppose a continuous non-negative function $g(x), ~x\in \R^I$, is fixed.
For $ x\in \cx^{(n)}$  define
\beql{eq-W-def}
G^{(n)}(x) = \int_0^{\infty} g(y^{(n)}(t)) dt, ~~ y^{(n)}(0)=x,
\end{equation}
where $y^{(n)}(\cdot)$ is the DFL starting from $x$.

\subsection{A tightness criterion}

Denote by $\nabla_z G^{(n)}(x)$ the directional derivative of $G^{(n)}$
at $x\in \cx^{(n)}$ in the direction of vector $z\in \R^I$:
$$
\nabla_z G^{(n)}(x) \doteq
\lim_{\delta\downarrow 0} \frac{1}{\delta} [G^{(n)}(x+z\delta)-G^{(n)}(x)],
$$
when the limit exists. (To be precise, if $x$ in on the boundary of $\cx^{(n)}$,
it is also required that the direction $z$ from $x$ points into $\cx^{(n)}$.)
Then, $\nabla_{z_*} [\nabla_z G^{(n)}](x)$ is the
second derivative, first in the direction $z$ and then $z_*$.

\begin{thm}
\label{th-criterion}
Suppose that for any $C_1 >0$,
there exists a function $g(x), ~x\in \R^I$, and a constant $C_2>0$,
such that the following conditions (i) and (ii) hold.

(i) Function $g(x)$  is  Lipschitz continuous,
non-negative and such that $g(x) \to \infty$ as $x\to\infty$.

(ii) For any $n$, the function $G^{(n)}(x), ~x\in\cx^{(n)}$, is finite for all $x$, and 
it has continuous gradient $\nabla G^{(n)}(x)$; for any $n$,
any $x$ and any fixed unit-length vectors $z,z_*\in \R^I$,
\beql{eq-second-deriv-bound}
\limsup_{\delta\downarrow 0} \frac{1}{\delta} 
\left| \nabla_z G^{(n)}(x+z_* \delta) -  \nabla_z G^{(n)}(x)\right|
\le C_1 g(x) + C_2;
\end{equation}
\beql{eq-W-to-inf}
 G^{(n)}(x) \to \infty, ~~x\to\infty, ~~~~\mbox{uniformly in $n$}.
\end{equation}
Then, 
$$
\sup_n \E  g(\hat X^{(n)}(\infty)) < \infty.
$$
\end{thm}

The second derivative condition \eqn{eq-second-deriv-bound} is the key one.
It implies that $|\nabla_{z_*} [\nabla_z G^{(n)}](x)| \le C_1 g(x) + C_2$ if this second derivative exists.
An equivalent form of  \eqn{eq-second-deriv-bound} is as follows:
for any compact set $D\subseteq \cx^{(n)}$ and any
unit-length vector $z\in \R^I$, the first derivative $\nabla_z G^{(n)}$ within $D$ is Lipschitz continuous
with constant
$$
C_1 \max_D g(x) + C_2.
$$

{\em Proof of Theorem~\ref{th-criterion}.} Let us fix a constant $C_1 >0$, 
and then the corresponding function $g$ and constant 
$C_2$, so that (i) and (ii) hold. (We will specify the choice of $C_1$ below.)

By definition of $G^{(n)}$ and its assumed continuous differentiability,
\beql{eq-W-grad}
\nabla G^{(n)}(x) \cdot v^{(n)}(x) =  \nabla_{v^{(n)}(x)} G^{(n)}(x)   =  - g(x).
\end{equation}

Let $A^{(n)}$ denote the (infinitesimal) generator of the Markov process  $X^{(n)}$. 
(See e.g. \cite{Ethier_Kurtz}, Sections 4.1 and 1.1, for the definitions of the operator semigroup and corresponding
generator of a Markov process. In our case, the semigroup is defined on the Banach space of bounded functions
$h(x), x \in \cs^{(n)},$ with norm $\sup_x |h(x)|$.)
For any fixed $k>0$, the function $G^{(n),k} \doteq G^{(n)} \wedge k$ is such that
it has constant value $k$ for all states $x\in  \cs^{(n)}$ except a finite subset $\cs^{(n),k}$ (where the value is less than $k$).
Therefore, there is only a finite number of possible
state transitions that may change the value of $G^{(n),k} (\hat X^{(n)}(t))$, namely
the transitions to or from the states in $\cs^{(n),k}$; the rates of such transitions are obviously uniformly bounded.
Using this property,
it is easy to verify directly that function $G^{(n),k}$  is within the domain of $A^{(n)}$, that is
$$
\lim_{t\downarrow 0} (1/t) E_x [G^{(n),k} (\hat X^{(n)}(t)) - G^{(n),k} (x)] =
A^{(n)} G^{(n),k} (x) = \sum_{x'} [G^{(n),k} (x')-G^{(n),k} (x)] \nu^{(n)}(x'-x),
$$
where $E_x$ denotes the expectation conditioned on $\hat X^{(n)}(0)=x$, and the limit is
uniform in $x$.
This in turn implies
\beql{eq-steady-state}
\E A^{(n)} G^{(n),k} (\hat X^{(n)}(\infty)) = 0.
\end{equation}
(See also \cite{Gamarnik_Stolyar}, page 31, for this property and argument in a very similar setting.)
For any $x\in \cs^{(n),k}$ we have
$$
A^{(n)} G^{(n),k}(x) \le \nabla G^{(n)}(x) \cdot v^{(n)}(x) + r^{(n)}(x)(1/2) h^{(n)}(x) (R_2/\sqrt{n})^2,
$$ 
where $R_2/\sqrt{n}$ is the maximum possible size of one jump of  the process, $r^{(n)}(x) \le R_1 n$ is the total transition rate from state $x$,
and the second-term coefficient $h^{(n)}(x)$ is bounded as $|h^{(n)}(x)| \le C_1 [g(x)+\kappa_1] +  C_2 = C_1 g(x) +C_1 \kappa_1 + C_2< \infty$.
(The constant $\kappa_1$ appears here, because we need an upper bound on the second derivative
in the $R_2$-neighborhood of point $x$, and we use the fact that $g(\cdot)$ is Lipschitz.)
Recalling also \eqn{eq-W-grad}, we obtain
$$
A^{(n)} G^{(n),k}(x) \le -g(x) + (1/2) R_1 R_2^2 [C_1 g(x) + C_1 \kappa_1 + C_2].
$$ 
We now specify the choice of $C_1$: it is sufficiently small so that (for any $x\in \cs^{(n),k}$)
$$
A^{(n)} G^{(n),k}(x) \le  -\epsilon g(x) +\kappa_2, ~~\mbox{for some} ~ \epsilon>0, ~ \kappa_2>0.
$$
(Function $g(\cdot)$ and constants $\kappa_1, C_2$ depend on $C_1$; therefore,
 constant $\kappa_2$ depends on the chosen $C_1$.). 
Obviously, if $x\in \cs^{(n)} \setminus \cs^{(n),k}$,
which is equivalent to $G^{(n),k}(x)=k$ and equivalent to $G^{(n)}(x) \ge k$, then
$$
A^{(n)} G^{(n),k}(x) \le 0.
$$
From these bounds and \eqn{eq-steady-state} we obtain
$$
\E [ -\epsilon g(\hat X^{(n)}(\infty)) +\kappa_2] I\{ G^{(n)}(\hat X^{(n)}(\infty)) < k\} \ge 
\E A^{(n)} G^{(n),k}(\hat X^{(n)}(\infty)) I\{ G^{(n)}(\hat X^{(n)}(\infty)) < k\} \ge
0,
$$
or
$$
\E  g(\hat X^{(n)}(\infty))  I\{ G^{(n)}(\hat X^{(n)}(\infty)) < k\} \le  \kappa_2/\epsilon.
$$
Letting $k\to \infty$, by monotone convergence,
$$
\E  g(\hat X^{(n)}(\infty)) \le  \kappa_2/\epsilon,
$$
where the constant in the RHS is independent of $n$.
$\Box$

\section{Proof of Theorem~\ref{th-tightnessN}}
\label{sec-proof-N}

We will prove Theorem~\ref{th-tightnessN} by choosing 
specific function $g(\cdot)$ and then verifying 
(in Theorem~\ref{thm-derivative-bound})
the assumptions of Theorem~\ref{th-criterion} for $N$-system.

In this section, we study properties of DFL trajectories and their $G^{(n)}$-functionals, for a system with a fixed 
scaling parameter $n$. We will drop upper index $(n)$ from now on. So, for example, will write simply $\cx$ and $y(t)$ instead of 
$\cx^{(n)}$ and $y^{(n)}(t)$, respectively. (However, the expressions may contain $n$ as a variable.) From this point on in the paper, we say that $C$ is a {\em universal
constant} if $C$ depends only on the system parameters 
$\lambda_i$, $\psi_{ij}$, $\mu_{ij}$, $b$,
but does {\em not} depend on scaling parameter $n$.
(If the sequence of systems is defined as in Remark~\ref{remark-rounding},
then a universal constant $C$ depends on the system parameters 
$\lambda_i$, $\psi_{ij}$, $\mu_{ij}$, $b$, but {\em not} on $n$ and {\em not} on the sequences $\psi_{ij}^{(n)}$
and $b^{(n)}$.)

\subsection{Basic DFL properties. First derivatives of DFLs and the Lyapunov function} 
\label{sec-basic-dfl}

In this subsection we first establish some basic properties of DFLs and their directional (Gateaux) derivatives.
Then we specify function $g(\cdot)$, and obtain the expressions for the first derivatives of the corresponding function $G(\cdot)$. 
(All results of this subsection hold for systems far more general than N-system.
In particular, they still hold for the systems under the
{\em Leaf Activity Priority} LAP discipline in \cite{SY2012,St2013_tightness}, in the Halfin-Whitt regime; our priority discipline for the $N$-system is a special case of LAP.)

The DFL trajectories $y(\cdot)$ have the following structure. 
Recall that $v(x)$ is (uniformly in $n$) Lipschitz continuous on the entire $\cx$.
There is a finite number $M$ (same for any $n$) of domains,  
indexed by $m=0,\ldots,M-1$; within each of them $v(x)$ is a given linear function. More precisely, 
the DFL satisfies a linear ODE
$$
(d/dt) y(t) = v(y(t)) = u^m y(t) + a^m,
$$
where $u^m$ is a constant $I\times I$ matrix (same for each $n$), and $a^m$ is a constant vector (depending on $n$).
Informally speaking, a domain is determined by which service pools a fully occupied and/or which queues are non-empty.

Formally, the domains are easier to define (and think of)
in terms of unscaled quantities $X_1\ge 0$ and $X_2\ge 0$, and unscaled pool sizes $B_1 = \psi_{11} n$
and $B_2 = \psi_{12}n + \psi_{22}n + b\sqrt{n}$.
Each domain is defined by a combination of the directions of three strict inequalities:
\beql{eq-ineqX1}
X_1 < B_1 ~~\mbox{or}~~ X_1 > B_1,
\end{equation}
\beql{eq-ineqX2}
X_2 < B_2 ~~\mbox{or}~~ X_2 > B_2,
\end{equation}
\beql{eq-ineqX1X2}
X_1 + X_2 < B_1+ B_2 ~~\mbox{or~~} X_1 + X_2 > B_1+ B_2.
\end{equation}
However, we exclude two combinations, or conditions, $(X_1 < B_1, ~X_2 < B_2, ~ X_1 + X_2 > B_1+ B_2)$
and $(X_1 > B_1, ~X_2 > B_2, ~ X_1 + X_2 < B_1+ B_2)$, because they produce the empty set;
and we replace ("merge") the conditions $(X_1 < B_1, ~X_2 > B_2, ~ X_1 + X_2 > B_1+ B_2)$
and  $(X_1 < B_1, ~X_2 > B_2, ~ X_1 + X_2 < B_1+ B_2)$ into one condition
$(X_1 < B_1, ~X_2 > B_2)$ because this condition alone determines the form of $v(x)$.
So, there are $M=5$ domains in total. 
The diffusion-scaling mapping $L_n$, defined by \eqn{eq-diff-scale},
 transforms them into 5 (diffusion-scale) 
domains, denoted $\cx^0, \ldots, \cx^4$. Note that the domains are defined
by strict inequalities, so they do not cover the entire space $\cx$. The domain closures
are $\bar\cx^1, \ldots, \bar\cx^5$, these do cover the entire $\cx$.
By these definitions, if a point belongs to the intersection of the closures of more than one domain,
then necessarily at least one of the equalities (in terms of unscaled quantities),
$X_1 = B_1$, $X_2 = B_2$, $X_1 + X_2 = B_1+ B_2$, holds.

In particular, consider the unscaled domain 
$(X_1 > B_1, ~X_2 < B_2, ~ X_1 + X_2 < B_1+ B_2)=(X_1 > B_1, ~ X_1 + X_2 < B_1+ B_2)$;
it is such that there are  no queues and pool 1 fully occupied. 
The corresponding diffusion-scaled domain is
$\cx^0=\{x\in \cx ~|~ x_1 > - \psi_{12}\sqrt{n},~ x_1+x_2 < b \}$.
In this domain
$v(x)=(- \mu_{12} x_1, - \mu_{22} x_2)$, i.e. $u^0=diag(- \mu_{12}, - \mu_{22})$ and $a^0=0$,
and therefore the components $y_1$ and $y_2$ evolve independently.
Moreover, there exists a universal constant $\alpha>0$, 
such that if $y(t)$ starts from a point
$y(0) \in \cx^{0,\alpha} \doteq \{\|x\| \le \alpha\} \subset \cx^0$, then $y(t)$ never leaves domain $\cx_0$,
which in turn means that the trajectory is simply
$$
y_i(t) = y_i(0) e^{-\mu_{i2} t}, ~~~i=1,2.
$$
From now such constant $\alpha$ and the corresponding sub-domain $\cx^{0,\alpha}$
will be fixed.

Consider one more unscaled domain 
$(X_1 > B_1, ~X_2 < B_2, ~ X_1 + X_2 > B_1+ B_2)$.
Here, pool 1 is fully occupied by type 1, pool 2 is fully occupied
by $X_2$ type 2 customers and $B_2 - X_2$ type 1 customers,
and $X_1 - B_1 -(B_2 - X_2)= X_1 + X_2 - B_1 -B_2 > 0$ type 1 customers waiting in the 
queue.
On the diffusion scale, the domain (let us label it $m=1$) is:
$\cx^1 = \{x\in \cx ~|~ x_1 > -\psi_{12}\sqrt{n}, ~
x_2 <  \psi_{12} \sqrt{n} + b, ~ x_1+x_2 > b\}$,
and we have
$$
v(x)=(\ (-b+x_2)\mu_{12},~ - \mu_{22} x_2),
$$
with the corresponding $u^1$ and $a^1$.
For the remaining 3 domains the $v(x)$ is determined similarly.

The equations for a
DFL $y(\cdot)$ can be summarized as follows. The trajectory of $y_2$ is 
not affected by $y_1$ and satisfies ODE
\beql{eq-y2}
(d/dt) y_2 = -\mu_{22} [y_2 \wedge (\psi_{12}\sqrt{n} + b)].
\end{equation}
If $y_1 \le - \psi_{12}\sqrt{n}$ (which corresponds to unscaled condition $X_1 \le B_1$),
\beql{eq-y11}
(d/dt) y_1 = -\mu_{11} (y_1 + \psi_{12}\sqrt{n}) + \mu_{12} \psi_{12}\sqrt{n} \ge \mu_{12} \psi_{12}\sqrt{n}.
\end{equation}
If $y_1 \ge - \psi_{12}\sqrt{n}$ ( $X_1 \ge B_1$) and 
$y_2 \ge \psi_{12}\sqrt{n} +b $ ( $X_2 \ge B_2$),
\beql{eq-y12}
(d/dt) y_1 =  \mu_{12} \psi_{12}\sqrt{n}.
\end{equation}
If $y_1 \ge - \psi_{12}\sqrt{n}$ ( $X_1 \ge B_1$), 
$y_2 \le \psi_{12}\sqrt{n} +b $ ( $X_2 \le B_2$), and $y_1+y_2\le b$ ($X_1+X_2\le B_1+B_2$),
that is in domain $\bar \cx^0$,
\beql{eq-y13}
(d/dt) y_1 =  -y_1 \mu_{12}.
\end{equation}
If $y_1 \ge - \psi_{12}\sqrt{n}$ ( $X_1 \ge B_1$), 
$y_2 \le \psi_{12}\sqrt{n} +b $ ( $X_2 \le B_2$), and $y_1+y_2\ge b$ ($X_1+X_2\ge B_1+B_2$),
that is in domain $\bar \cx^1$,
\beql{eq-y14}
(d/dt) y_1 =  (-b+y_2) \mu_{12}.
\end{equation}

For a given fluid trajectory, let us call time point $t\ge 0$ a {\em switching point}
if $y(t)$ belongs to the intersection of two or more closed domains $\bar \cx^m$.
(i.e. it is on a boundary separating different domains).

\begin{lem}
\label{lem-basic1}
For some universal constants $T>0$, $C'>0$ and (integer) $K'>0$, DFL trajectories $y(\cdot)$ 
satisfy the following conditions. [For a DFL $y(\cdot)$, $x=y(0)\in \cx$ denotes its initial state.]\\
(i) Let $\tau\ge 0$ be the first time a DFL reaches set $\cx^{0,\alpha}$. Then,  $\tau\le T \|x\|$.
(This, in particular, means that $y(t) \to 0, ~t\to\infty$, and, moreover, the convergence is exponentially fast.)
In addition, $\|y(\cdot)\| \le C' \|x\|$. \\
(ii) DFL $y(\cdot)$ depends on its initial state $x$ continuously, in the sense of $\|y(\cdot)\|$-norm.\\
(iii) DFL $y(\cdot)$ has at most $K'$ switching points, $t_1 < t_2 < \ldots < t_K$, $ 0 \le K\le K'$,
and $t_K < \|x\|T$. 
 Moreover,
the set of switching points is upper semicontinuous in $x$; namely,
as $x \to x^*$, the limiting points of the set 
of switching points are within the set of switching points for initial state $x^*$.\\
(iv) For any interval $[C_3,C_4]$, not containing $0$, there exists a constant $T_3>0$ (independent of $n$), 
such that the total time the condition $y_i(t) \in [C_3,C_4]$ holds for at least one $i$, is upper bounded by $T_3$.
\end{lem}

{\em Proof of Lemma~\ref{lem-basic1}.} 
Given equation \eqn{eq-y2},
condition $y_2(t) = \psi_{12}\sqrt{n} + b$ ( $X_2=B_2$) can hold at most at one point
$t_2 \ge 0$, which will be a switching point. 
Similarly, by \eqn{eq-y11}, there is at most one point $t_1 \ge 0$, at which condition $y_1(t) = - \psi_{12}\sqrt{n}$
(corresponding to $X_1=B_1$) can hold, and if so, it will be a switching point. 

Denote $t'=t_1 \vee t_2$.
It is easy to see that for some universal constant $\kappa>0$,
\beql{eq-tprime-bound}
t' \le  \kappa\|x\|, ~~\|y(t')\| \le \kappa \|x\|.
\end{equation}
Indeed,  $|y_2(t)|$ is non-increasing in $[0,\infty)$, and $t_2 \le |x_2|/[(\psi_{12}\sqrt{n} +b) \mu_{22}] \le | x_2| /  [\psi_{12} \mu_{22} \sqrt{n}]$. In the interval $[0,t_1]$, $y_1(t)$ is negative  non-decreasing, and then $|y_1(t)|$ is non-increasing;
and $t_1 \le  |x_1|/[\psi_{12}  \mu_{12} \sqrt{n}]$. If $t_2 > t_1$, then in the interval $[t_1,t_2]$, 
$(d/dt) y_1(t) = \psi_{12}  \mu_{12} \sqrt{n}$, and therefore $|y_1(t_2)-y_1(t_1)| \le \psi_{12}  \mu_{12} \sqrt{n} t_2$;
given the bound on $t_2$, we see that $|y_1(t_2)-y_1(t_1)|$ is upper bounded by $|x_2|$ times a universal constant.
These observations imply \eqn{eq-tprime-bound}.

For all $t>t'$, conditions $y_2(t) < \psi_{12}\sqrt{n} + b$ ( $X_2<B_2$) 
and $y_1(t) > - \psi_{12}\sqrt{n}$ ( $X_1>B_1$) hold. Therefore, $y(t)$ can be only in one of the
two domains $\bar \cx^0$ or $\bar \cx^1$, depending on whether $y_1 + y_2 \le b$ (no queues) or $y_1 + y_2 \ge b$ (queue size $y_1 + y_2 - b$ of type 1).
It is easy to see from equations $(d/dt)y_2 = -\mu_{22} y_2$, \eqn{eq-y13}, \eqn{eq-y14}, 
that if $y(t)$ is in $\bar \cx^1$, then the trajectory eventually leaves $\bar \cx^1$
and can never return.
This implies that at most two transitions between $\cx^0$ and $\cx^1$ can occur after $t'$.
Specifically, either the trajectory stays in $\cx^0$, or it is in $\cx^1$ and then $\cx^0$, or it is in 
$\cx^0$  then $\cx^1$  then $\cx^0$. The boundary cases are also possible; for example, the trajectory 
may stay in the open domain $\cx^0$ at all times, except at exactly one point $t\ge t'$ it "touches" the boundary,
i.e.  $y_1 + y_2 = b$. To summarize, after $t'$ there are at most two switching points.

Denote by $t''$ the first time $t\ge t'$ when $\|y_2(t)\| \le \alpha/4$.
We have $t''-t' = 0\vee (1/\mu_{22}) \log [\|y_2(t')\|/(\alpha/4)] \le \kappa_1 \|x\| + \kappa_2$,
for some universal $\kappa_1$ and $\kappa_2$. ($\kappa_2$ depends on $\alpha$, which in turn is universal.)
In the interval $[t',t'']$ the value of $|y_1|$ cannot increase by more than
$\kappa_3 |y_2(t')| \le \kappa_4 \|x\|$, for universal $\kappa_3, \kappa_4>0$.
(If $y_1\le 0$, then $(d/dt)y_1 \ge 0$. If $y_1\le 0$, then $(d/dt)y_1 \le \mu_{12} |y_2|$,
and recall that $(d/dt)y_2 = - \mu_{22} y_2$.)
Therefore, $|y_1(t'')| \le \kappa_5 \|x\|$, for universal $\kappa_5>0$.
Starting $t''$, if type 1 has non-zero queue, $(d/dt)|y_1| = (d/dt)y_1 \le - \kappa_6 < 0$,
for universal $\kappa_6>0$;
and if type 1 does not have queue, then $(d/dt)|y_1| = - \mu_{12} |y_1|$.
Consider the first time $t'''\ge t''$ when $|y_1| \le \alpha/4$. 
We conclude that $t''' \le T \|x\| + \kappa_7$
and $\sup_{[0,t''']} \|y(t)\| \le C' \|x\|$ 
for some universal positive constants $T,C',\kappa_7$.
Obviously, $t'''\ge \tau$, so that
$\tau \le T \|x\| + \kappa_7$.
However, if $\|x\| \le \alpha$, i.e. $y(0)=x$ is already in $\cx^{0,\alpha}$, then obviously $\tau=0$.
Therefore, in the bound $\tau \le T \|x\| + \kappa_7$, we can drop $\kappa_7$ by rechoosing $T$, if necessary.

For future reference, we also make the following observation. 
{\em Suppose, $\mu_{12}=\mu_{22}$.
Then, there can be at most one switching point after time $t'$, let us call it $t_3\ge t'$, and it is such that
$y(t)\in \cx^0$ for all $t>t_3$.}
Indeed, in this case, in the domain $\bar \cx^0$, we have simply 
$(d/dt)[y_1+y_2] = - \mu_{22}[y_1+y_2]$.

Let us prove properties (i)-(iv). In fact, (i) has been proved already. For a given $x$, let us choose $\tau'$ such that $\tau < \tau'$
for all initial states sufficiently close to $x$. (On a finite interval $[0,\tau']$, $y(\cdot)$ depends on the initial state continuously,
because it is a solution to an ODE with Lipschitz continuous RHS.) But, for $t\ge \tau'$, the DFL with any initial state close
to $x$ is such that $y(t)\in \cx^{0,\alpha}$; this implies uniform convergence across all $t\ge 0$, which proves (ii).
The part of property (iii), stating that there is at most $K'$ switching points, all of which are smaller than $\tau \le T\|x\|$, has already been
proved, in fact we specified that $K'\le 4$. Then, the upper semicontinuity of the set of switching points follows from continuity of
trajectories w.r.t. initial state; this proves (iii). 
Consider a fixed interval $[C_3,C_4]$, not containing $0$. 
It is clear from \eqn{eq-y2} that 
$y_2(t)$ can spend only a finite time within $[C_3,C_4]$.
Now, $y_1(t)$ can be in $[C_3,C_4]$ only after time $t_1$,
and then in every domain the trajectory visits $y_1(t)$
satisfies one of the equations \eqn{eq-y12}-\eqn{eq-y14}. 
If we examine each of these equations (and recall that \eqn{eq-y14} holds
within domain $\bar \cx^1$, where $(d/dt) y_2 = -\mu_{22} y_2$),
we see that {\em even if the equation were to hold up to infinite time}, 
$y_1(t)$ can spend only a finite time within $[C_3,C_4]$. 
And there is only a finite, uniformly bounded
number of domains that a trajectory can visit.
This proves (iv).
$\Box$

Next, let us consider the first-order 
dependence of DFL on the initial state.
Let $y(t;x)$ denote $y(t)$ with initial state $y(0)=x \in \cx$. 
For any $x\in \cx$ and any direction $z\in \R^I$ (which does not point outside $\cx$),
we use the following notation for the directional (Gateaux) derivative of $y(t;x)$ at $x$ in the direction $z$:
$$
\nabla_z y(\cdot;x) \doteq \lim_{\delta\downarrow 0} \frac{1}{\delta} [y(\cdot;x+ z\delta) - y(\cdot;x)].
$$

\begin{thm}
\label{thm-variation}
(i) For any fixed $x \in \cx$ and a fixed vector $z$,
the directional derivative
$$
\xi(\cdot) = \xi(\cdot;x,z) = \nabla_z y(\cdot;x) 
$$
exists. It has the following structure. Let $0 < t_1 < t_2 < \ldots < t_K$ be the switching points of $y(\cdot;x)$.
Then, $\xi(0)=z$, and in each interval $[0,t_1], [t_1,t_2], \ldots, [t_K,\infty)$, $\xi$ satisfies linear homogeneous ODE
$$
(d/dt) \xi = u^m \xi,
$$
where matrix $u^m$ is the matrix $u$ for the domain $\bar \cx^m$ containing $y(t;x)$.\\
Solutions $q(t), ~t\ge 0,$ to the equation $(d/dt)q = u^m q$, for any $m$,
are such that 
\beql{eq-factor}
\|q(\cdot)\| \le C_5 \|q(0)\|
\end{equation}
for a universal constant $C_5>0$.

(ii) The derivative $\xi(\cdot;x,z)$ depends on $(x,z)$ continuously.

(iii) There exists a universal constant $C_6>0$,
such that 
$$
\|\xi(\cdot;x,z)\| \le C_6 \|\xi(0;x,z)\|= C_6 \|z\|.
$$
\end{thm}

{\em Proof.}
The proof of (i) relies on the following observations.\\
 (a) In any time interval, where both 
$y(t;x+ z\delta)$ and $y(t;x)$ are within same domain $\bar \cx^m$, 
they are governed by the same  ODE $(d/dt) y = v^m(y)$, and therefore their difference
 $\Delta y(t) = y(t;x+ z\delta)-y(t;x)$, is governed by the linear homogeneous ODE
$(d/dt) \Delta y = u^m \Delta y$. Moreover, it is easy to check that 
within any domain  $\bar \cx^m$ the corresponding matrix $u^m$ is such that
$\|\Delta y(t)\|$  can increase at most by some universal factor $C_8$. 
Indeed, consider $\Delta y_2$ first, and then $\Delta y_1$.
The equation for $\Delta y_2$ is either 
\beql{eq-delta-der2a}
(d/dt) \Delta y_2 = -\mu_{22} \Delta y_2
\end{equation}
or 
\beql{eq-delta-der2b}
(d/dt) \Delta y_2 = 0;
\end{equation}
in either case $|\Delta y_2|$ cannot increase.
The equation for $\Delta y_1$ is
$$
(d/dt) \Delta y_1 = u^m_{11} \Delta y_1 + u^m_{12} \Delta y_2,
$$
where $u^m_{11}=0$ or $u^m_{11}=-\mu_{11}$ or $u^m_{11}=-\mu_{12}$;
we also note that if $\Delta y_2$ satisfies \eqn{eq-delta-der2b}
then necessarily $u^m_{12} = 0$. We see that in any case, in any time interval,
$|\Delta y_2(t)|$ is upper bounded by the initial $\|\Delta y_2\|$ 
times a universal constant. This observation, in particular, 
proves \eqn{eq-factor}.\\
(b) The total length of "switching intervals", where $y(t;x+ z\delta)$ and $y(t;x)$ belong to different domains vanishes
as $\delta\to 0$ (by upper semicontinuity of the set of switching points),
and therefore the total change of $\Delta y(t)$ within those intervals is "small".
More precisely,  let $t$ be fixed and 
$[\theta_1,\theta_2]$ be a switching interval such that
$\theta_1, \theta_2 \to t$. Then, 
$\|\Delta y(\theta_2) - \Delta y(\theta_1)\|/\|\Delta y(\theta_1)\| \to 0$,
because $v(x)$ is Lipschitz.\\
Combining observations (a) and (b), and further observing that the number of 
intervals where both 
$y(t;x+ z\delta)$ and $y(t;x)$ are within same domain $\bar \cx^m$
(i.e. outside the switching intervals) is upper bounded, we take the
$\delta\downarrow 0$ limit to obtain (i).

(ii) This follows from the upper semicontinuity of the set of switching points 
on $x$.

(iii) By \eqn{eq-factor},  in any domain $\|\xi(t)\|$ can increase at
most by some factor $C_5$. There is only a finite number of domains that $y(t)$ visits. This proves (iii).
$\Box$

We now introduce a specific function $g$, which we will use in the definition 
\eqn{eq-W-def}
of the Lyapunov function.

\begin{definition}
\label{def2}
Let parameter $C > 0$ be fixed. Let a function $f(\eta)$ of real $\eta$ be fixed, which satisfies the following conditions.
It is a non-negative, even, convex, twice continuously differentiable, 
 $f(\eta)=0$ for $\eta\in [-C,C]$, $f'(\eta) = -1$ for $\eta\le -C-1$,
$f'(\eta) = 1$ for $\eta\ge C+1$. 
(Such a function can be defined explicitly. 
Since $C$ is a parameter, essentially, we fix the shape of function $f(C+ \zeta), ~\zeta\ge 0$.)
Note that both $f'$ and $f''$ are uniformly bounded, and $f''=0$ outside
of the intervals $[-C-1,-C]$ and $[C,C+1]$.
Then, let
$$
g(x) = \sum_i f(x_i).
$$
Obviously, $\left| f(\eta)-|\eta|~ \right|$ is uniformly bounded by a constant, and then so is $\left| g(x)-\|x\| ~ \right|$.
\end{definition}

Then, by \eqn{eq-W-def} we have $G(x)=\sum_i G_i(x)$, where 
\beql{eq-Wi-def}
G_i(x) = \int_0^{\infty} f(y_i(t)) dt, ~~ y(0)=x.
\end{equation}
Clearly, $G(x)$ is finite for any $x$, because (by Lemma~\ref{lem-basic1}(i)) $y(t)\to 0$ and therefore
(by definition of $f(\cdot)$) $f(y_i(t))=0$ for all large $t$.

\begin{thm}
\label{thm-derivatives}
For each $i$ the following holds.
For any $x\in \cx$ and any direction vector $z$,
\beql{eq-deriv-1}
\nabla_z G_i(x) = \int_0^\infty f'(y_i(t;x)) \xi_i(t;x,z) dt.
\end{equation}
Function $\nabla_z G_i(x)$ is continuous
in $(x,z)$.
\end{thm}

{\em Proof.} 
Expression \eqn{eq-deriv-1} follows from Theorem~\ref{thm-variation}(i)
and the fact that $f'$ is continuous bounded. The continuity 
of $\nabla_z G_i(x)$ is obtained using 
Theorem~\ref{thm-variation}(ii) and Lemma~\ref{lem-basic1}(i,ii).
$\Box$

\subsection{Second derivative bounds for the Lyapunov function}
\label{subseq-second-deriv}

\begin{thm}
\label{thm-derivative-bound}
The assumptions of Theorem~\ref{th-criterion} hold. Specifically, for any $C_1>0$, there exist constants $C>0$ and $C_2>0$
such that conditions (i) and (ii) in Theorem~\ref{th-criterion} hold for the 
function $g$ in Definition~\ref{def2} with parameter $C$.
\end{thm}

Note that for a function $g$ satisfying Definition~\ref{def2},
 condition (i) of Theorem~\ref{th-criterion} holds automatically.
Condition \eqn{eq-W-to-inf} is also automatic given the definition of $G$ and basic properties of DFL,
namely the fact that the time for a DFL to reach a given compact set increases to infinity as $x\to\infty$.
Therefore, to prove Theorem~\ref{thm-derivative-bound}, it remains to prove
condition \eqn{eq-second-deriv-bound}, and it suffices to prove it separately for $G_i$, $i=1,2$ (see \eqn{eq-Wi-def}).
We will do this first for the case $\mu_{22} \ne \mu_{12}$, and then for $\mu_{22} = \mu_{12}$.
(The proof of condition \eqn{eq-second-deriv-bound} in this section
applies to the  N-system, as well as its generalization described
in Section~\ref{sec-generalization}. It does {\em not} apply for LAP discipline.)

For a given $x$ and a time $\tau^*>0$, denote by $S(\tau^*;x)$ the set of time points,
consisting of $\tau^*$ and all switching points $0 \le t<\tau^*$ of the DFL
$y(\cdot;x)$.

\begin{lem}
\label{lem-switchover}
Suppose  $\mu_{22} \ne \mu_{12}$.
For any $\epsilon>0$ there exists a sufficiently large $C_7>0$,
such that, for all sufficiently large $n$,
 the following holds for any fixed $x$ and any unit-length vector $z$.
Let $\tau_7$ be the first time the DFL $y(\cdot;x)$
hits set $\{\|y\| \le C_7\}$. Then
for all sufficiently small $\delta>0$, any point in $S(\tau_7;x+z\delta)$
is within distance at most $\epsilon \delta$
from a point in $S(\tau_7;x)$.
\end{lem}

{\em Proof.} 
Consider a switching point  $t\in S(\tau_7;x)$ of 
DFL $y(\cdot)=y(\cdot;x)$. 
By definition of $\tau_7$, it
is such that $\|y(t)\| \ge C_7$.
The switching point is on the boundary of multiple domain closures,
and therefore one or more equalities
\beql{eq-dom-boundary}
y_1(t)=-\psi_{12}\sqrt{n}, ~~
y_2(t)=\psi_{12} \sqrt{n} + b, ~~ y_1(t)+y_2(t) = b,
\end{equation}
defining the domain boundaries, hold.
If the first or second equality holds,
then $|y'_i(t)|$ is large for large $n$.
If $y_1(t)+y_2(t) = b$, then for $t$ to be a switching point, it is necessary that $y(t)\in \bar\cx^0$;
then $y'_1(t)+y'_2(t) = -\mu_{12} y_1(t) -\mu_{22} y_2(t)=  -(\mu_{12} -\mu_{22}) y_1(t) - \mu_{22}b $; conditions
 $y_1(t)+y_2(t) = b$ and $\|y(t)\|=|y_1(t)| + |y_2(t)| \ge C_7$ imply that if $C_7$ is large then so is $|y_1(t)|$,
and then $|y'_1(t)+y'_2(t)|$ is large as well. We conclude that if any of the three equalities 
\eqn{eq-dom-boundary} holds, then for all $n\ge  n'$
we have $|y'_1(t)| \ge \kappa_1$ or  $|y'_2(t)| \ge \kappa_1$ or  $|y'_1(t)+y'_2(t)| \ge \kappa_1$,
respectively, where the constant $\kappa_1>0$ can be made arbitrarily large by choosing large enough $n'$ and $C_7$.
This means that, first, the domains in which the trajectory $y(\cdot;x)$ is in before and after
the switching point $t$ are uniquely defined. Second, since 
the distance between $y(\cdot;x+z\delta)$ and $y(\cdot;x)$ does not exceed
$\kappa_2 \delta$ at all times, where $\kappa_2 >0$ is a universal constant
(this follows from Theorem~\ref{thm-variation}),
and $v(\cdot)$ is Lipschitz,
any point in $S(\tau_7;x+z\delta)$ must be within $2 \kappa_2\delta/\kappa_1$
of a point in $S(\tau_7;x)$. Since $\kappa_2$ is universal and
$\kappa_1$ can be made arbitrarily large (by choosing $C_7$ large), the result follows.
$\Box$

Recall that to prove Theorem~\ref{thm-derivative-bound}, it remains to prove
the second derivative condition \eqn{eq-second-deriv-bound}.
The proper second derivative may not exist, hence we must ``settle'' for the 
estimate \eqn{eq-second-deriv-bound}.
But, to illustrate the proof
that follows, let us write down the expression for the second derivative,
by formally applying $\nabla_{z_*}$ differentiation to \eqn{eq-deriv-1} 
(this expression is {\em not} used in the proof):
$$
\nabla_{z_*} \nabla_z G_i(x) =
$$
\beql{eq-deriv-2-formaA}
\int_0^\infty f''(y_i(t;x)) \xi_i(t;x,z_*) \xi_i(t;x,z)dt +
\end{equation}
\beql{eq-deriv-2-formaB}
\int_0^\infty f'(y_i(t;x)) \nabla_{z_*}\xi_i(t;x,z)dt.
\end{equation}

{\em Proof of Theorem~\ref{thm-derivative-bound}, case $\mu_{22} \ne \mu_{12}$.}
We choose small $\epsilon>0$ and then $C_7>0$ as in Lemma~\ref{lem-switchover}.
Then choose parameter $C>0$ of 
function $g$ large enough so that
any DFL starting from the set $\{\|y\| \le 2 C_7\}$
never hits set $\{\|y\| \ge C\}$. (We can do this by Lemma~\ref{lem-basic1}(i).)

For $i=1,2$ consider 
$$
\frac{1}{\delta} [f'(y_i(t;x+z_* \delta)) \xi_i(t;x+z_* \delta, z)
- f'(y_i(t;x)) \xi_i(t;x,z)]
$$
\beql{eq-key1}
= \frac{1}{\delta} [f'(y_i(t;x+z_* \delta)) - f'(y_i(t;x))] \xi_i(t;x,z) +
\end{equation}
\beql{eq-key2}
\frac{1}{\delta} f'(y_i(t;x+z_* \delta))[\xi_i(t;x+z_* \delta, z) - \xi_i(t;x, z)].
\end{equation}
(The integrals of the terms \eqn{eq-key1} and \eqn{eq-key2}, correspond to
the integrals \eqn{eq-deriv-2-formaA} and \eqn{eq-deriv-2-formaB}, respectively,
in the formal second derivative expression.)

Since $f(\cdot)$ has bounded second derivative, the term \eqn{eq-key1} converges (uniformly in $t$) to 
$$
f''(y_i(t;x)) \xi_i(t;x,z_*) \xi_i(t;x,z).
$$
The integral of this over $t\in [0,\infty)$
is bounded because the total time any trajectory spends
in the set $\{C\le \|y_i\| \le C+1\}$ is uniformly bounded (by Lemma~\ref{lem-basic1}(iv).)

In the term \eqn{eq-key2}, $f'(y_i(t;x+z_* \delta))$ is uniformly bounded.
Let $\tau_7$ be the first time $y(t; x)$ hits set $\{\|y\| \le C_7\}$. 
We claim that, uniformly in $t\in [0,\tau_9]$,
\beql{eq-claim}
\limsup_{\delta\downarrow 0} \frac{1}{\delta} [\xi_i(t;x+z_* \delta,z) - \xi_i(t;x,z)]
\le \epsilon \kappa,
\end{equation}
where $\kappa>0$ is a universal constant.
Indeed, let $t_1\in S(\tau_7;x)$ be the first (smallest) switching point of trajectory 
$y(\cdot;x)$. To be concrete, let us assume $t_1>0$. (The case $t_1=0$ is treated analogously.)
For a given $\delta$, we define a {\em switching interval} 
$[\theta_1^*,\theta_1^{**}]$ associated with $t_1$ as follows: $\theta_1^*$ is the minimum
of $t_1$ and those switching points of $y(\cdot;x+z\delta)$ that are
within distance $\epsilon \delta$ from $t_1$;
similarly, $\theta_1^{**}$ is the maximum
of $t_1$ and those switching points of $y(\cdot;x+z\delta)$ that are
within distance $\epsilon \delta$ from $t_1$. Obviously, 
$\theta_1^{**}-\theta_1^* \le 2\epsilon\delta$.
In the interval $[0,\theta_1^*]$, $\xi(t;x+z_* \delta, z) = \xi(t;x,z)$,
because they are governed by the ODE with {\em same} matrix $u^m$.
Within the switching interval, the ODEs for $\xi(t;x+z_* \delta, z)$
and $\xi(t;x, z)$ may have a different matrix $u^m$, but there is only
a finite number of those matrices; therefore, in $[\theta_1^*,\theta_1^{**}]$,
$\|\xi(t;x+z_* \delta, z) - \xi(t;x, z)\|$
can increase at most by $\kappa_1\|\xi(\theta_1^*;x, z)\|\epsilon\delta$, 
where $\kappa_1$ is a universal constant.
We then consider the second switching point $t_2$ and the associated
switching interval $[\theta_2^*,\theta_2^{**}]$. Note that
between the first and second switching intervals,
both $\xi(t;x+z_* \delta, z)$
and $\xi(t;x, z)$ are again governed by the ODE with same matrix $u^m$;
therefore the difference $\xi(t;x+z_* \delta, z) - \xi(t;x,z)$
is governed by the same ODE, and therefore in the interval $[\theta_1^{**},\theta_2^*]$ 
the value of
$\|\xi(t;x+z_* \delta,z) - \xi(t;x,z)\|$ can increase at most by a factor given by a
 universal constant $\kappa_2>0$ (by \eqn{eq-factor}).
At the end of the switching interval $[\theta_2^*,\theta_2^{**}]$,
the first-order component of
 $\| \xi(t;x+z_* \delta, z) - \xi(t;x, z)\|$
is upper bounded by 
$$
\kappa_2 \kappa_1 \|\xi(\theta_1^*;x,z)\|\epsilon\delta +
\kappa_1 \|\xi(\theta_2^*;x,z)\|\epsilon\delta.
$$
We consider the third switching point, and so on. 
We see that
the first-order 
component of $\|\xi(t;x+z_* \delta,z) - \xi(t;x,z)\|$
will be upper bounded by $\kappa \epsilon \delta$, for a sufficiently large universal $\kappa$.
(There will be also higher order terms $\delta^{\ell}, ~\ell \ge 2,$ 
with uniformly bounded coefficients.)
This proves claim \eqn{eq-claim}.

By Lemma~\ref{lem-basic1}(i),
 $\tau_7 \le \kappa_3\|x\|$, for a universal constant $\kappa_3>0$.
(We can always choose $C_7 \ge \alpha$, and then $\kappa_3=T$.)
Then, $\tau_7 \le \kappa_3 (g(x)+\kappa_4)$, where $\kappa_4$ may depend on the parameter
$C>0$ of function $g$.

Now, for all sufficiently small $\delta$, the integral of the term \eqn{eq-key2},
$$
\int_0^\infty \frac{1}{\delta} f'(y_i(t;x+z_* \delta))[\xi_i(t;x+z_* \delta,z) - \xi_i(t;x,z)]dt =
$$
$$
\int_0^{\tau_9} \frac{1}{\delta} f'(y_i(t;x+z_* \delta))[\xi_i(t;x+z_* \delta,z) - \xi_i(t;x,z)]dt,
$$
because $f'(y_i(t;x+z_* \delta))=0$ for $t\ge T$.
The absolute value of the latter integral is upper bounded by
$$
\kappa_3 (g(x)+\kappa_4)\kappa \epsilon = \kappa \kappa\epsilon g(x)+\kappa_3 \kappa_4 \kappa \epsilon.
$$
The constants $\kappa_3$ and $\kappa$ are universal, 
while $\kappa_4$ depends on $C$, which depends on $C_7$,
which depends on $\epsilon$. It remains to choose $\epsilon$ small enough
so that $\kappa_3 \kappa\epsilon < C_1$. Then the value of 
$\kappa_3 \kappa_4 \kappa \epsilon$, plus the corresponding upper bound
on the integral of \eqn{eq-key1}, gives constant $C_2$.
$\Box$

{\em Proof of Theorem~\ref{thm-derivative-bound}, case $\mu_{22} = \mu_{12}$.}
This case is treated the same way as $\mu_{12}\ne\mu_{22}$, with the following modifications.
If there is {\em no} switching point $t\in S(\tau_7;x)$,  associated with 
equality $y_1(t)+y_2(t)=b$, then the proof is unchanged.
Suppose there is a switching point $t\in S(\tau_7;x)$, associated with 
equality $y_1(t)+y_2(t)=b$. Then, in the notation of the proof 
of Lemma~\ref{lem-basic1}, we must have $t \ge t'$, and by the observation
we made in that proof, 
$t$ is the last switching
point, and therefore it is the only switching point associated with
equality $y_1(t)+y_2(t)=b$.
Moreover,
all the properties we established in the  $\mu_{22} \ne \mu_{12}$ case proof,
still apply to all switching points before $t$.
After time $t$,
the process stays within the domain $\cx_0$, and therefore
$(d/dt)[y_1(t)+y_2(t)] = -\mu_{22} [y_1(t)+y_2(t)]$.
In particular, in a small neighborhood of time $t$, 
$(d/dt)[y_1(t)+y_2(t)] \le - (b/2) \mu_{22} <0$.
These facts imply that
the switching interval, corresponding to switching time $t$,
is such that its end points are within $\kappa_5 \delta$ from $t$, for some universal constant
$\kappa_5>0$. This means that 
the contribution of this last switching interval, as well as of the remaining
time interval up to the time $\tau_7$, into the integral of  \eqn{eq-key2},
is upper bounded by a universal constant $\kappa_6>0$.
$\Box$

\section{Generalization of the $N$-system}
\label{sec-generalization}

\begin{figure}[htbp]
\begin{center}
\includegraphics[width=140mm]{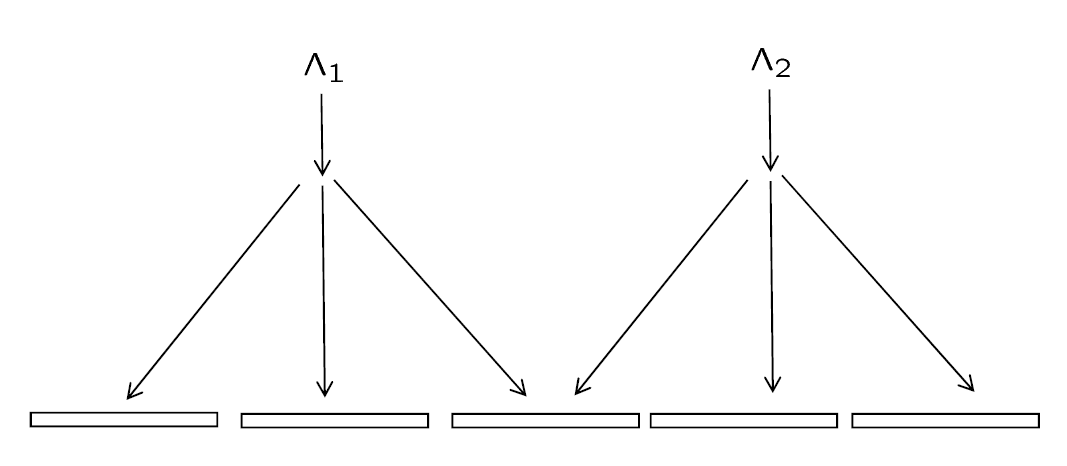}
\end{center}
\vspace{-0.2in}
\caption{A more general system.}
\label{fig:general-system}
\end{figure}

Theorem~\ref{th-tightnessN}, along with its proof, easily extend to the 
generalization of $N$-system, shown in Figure~\ref{fig:general-system},
in the Halfin-Whitt regime.
The system has two customer types and arbitrary number of server pools.
There is exactly one server pool that is flexible, i.e. can serve both types.
(On Figure~\ref{fig:general-system}, it is the pool in the middle.)
Each of the remaining pools is dedicated to service of either type 1 or 2.
(The two pools on the left in the figure are dedicated to type 1, while the
two pools on the right -- to type 2.) Each customer type has absolute preference for
its dedicated server pools, in some fixed priority order, over the flexible pool.
In the flexible pool, the absolute preemptive priority is given to one of the types.

The key features that the generalized system shares with the $N$-system are that
there are two customer types and only one flexible server pool, which can be shared
by the customers of different types. 
These features are exploited in Section~\ref{subseq-second-deriv}, where 
we estimated second derivatives of the Lyapunov function.
(We note again that all results in Section~\ref{sec-basic-dfl},
which concern with first derivatives,
hold for far more general systems, e.g. those under LAP discipline 
\cite{SY2012,St2013_tightness}.)
The behavior of the DFLs for the generalized system is more complicated,
simply because the number of state space domains can be very large.
However, as in the $N$-system, after a finite time all dedicated server pools
stay fully occupied, which means that the DFL dynamics depends only on
``what happens'' in the flexible pool. Consequently, our analysis goes through
with very minor adjustments.

\section{Discussion}
\label{discussion}

In this paper we address the problem of tightness of stationary distributions,
and the limit interchange, for flexible multi-pool service systems in the Halfin-Whitt 
regime. The behavior of such systems can be very complicated, which makes 
the problem challenging. This is, in particular, due to the difficulty of constructing
Lyapunov functions. Our approach uses a (family of) Lyapunov function(s),
defined as an integral functional of the drift-based fluid limits (DFL) $y(\cdot)$:
$G(x)=\int_0^\infty g(y(t)) dt$, $y(0)=x$.
The problem then reduces to studying the (first and second) derivatives
of a DFL -- and the corresponding integral $G(x)$ -- on the initial state $x$.
We apply this approach to show the tightness property for 
the $N$-model under a priority discipline. 

Both the approach 
and many parts of our analysis are quite generic and might be applicable
to other models as well. In this respect, note that there is a lot of flexibility
in choosing the ``distance'' function $g(\cdot)$. 
It might also be possible  to combine the approach with other techniques. 
For example, a Lyapunov function of the type we consider
could be defined and applied on a subspace, if it could be shown 
by other means that the stationary distributions concentrate (in appropriate sense) 
on that subspace. Exploring these directions may be a subject of future research.

\bibliographystyle{acmtrans-ims}
\bibliography{gen-lyap}

\appendix

\section{Positive recurrence proof}
\label{app1}

Let us drop superscript $(n)$. Consider the process with fixed initial state such that $X_2(0)=0$. Consider the sequence of time points 
$0 < t_1 < t_2 < \ldots$ at which $X_2(t)$ changes to $0$; let $t_0=0$. These are renewal points for $X_2(\cdot)$ viewed 
in isolation; $X_2(\cdot)$  is positive recurrent. The renewal interval durations $t_{\ell+1} - t_\ell$, $\ell=0,1,\ldots$, are of course i.i.d.
with some finite mean $T$. Let $A_\ell$ be the random number of type 1 arrivals into the system in the interval 
$(t_\ell, t_{\ell+1}]$; and $S_\ell$ be the random number of type 1 service completions  in the interval $(t_\ell, t_{\ell+1}]$,
{\em assuming that all servers (in both pools), not occupied by type 2 customers, serve type 1 customers}. 
Clearly, $(A_\ell, S_\ell)$ are i.i.d. across $\ell$, $\E A_\ell = \lambda_1 n T$, $\E S_\ell = \lambda_1 n T + b \mu_{12} \sqrt{n}T$,
$\E A_\ell - \E S_\ell = -  b \mu_{12} \sqrt{n}T < 0$.
Using these facts, it is easy to see that the discrete time Markov chain $X_1(t_\ell), \ell=0,1,2, \ldots$ is positive recurrent;
let $N<\infty$ denote the mean return time to $0$ for this chain. This implies that, for the original continuous-time
process $(X_1(t),X_2(t))$, the mean time to return to state $(0,0)$ is upper bounded by $NT$.
We omit further details.

\end{document}